\documentclass[11pt]{article}
\usepackage{amsmath,amssymb,mathrsfs}
\usepackage[french]{babel}
\usepackage[ansinew]{inputenc}
\usepackage{makeidx}
\usepackage{graphicx}
\usepackage{graphics}
\usepackage{newlfont}
\usepackage[T1]{fontenc}
\usepackage{geometry}
\usepackage{fancyhdr}
\makeindex  
 \advance\topmargin 0cm

\long\def\cor #1#2 {\par{\bf Corollary #1. }{\it #2}\par}
\long\def\lemme #1#2 {\par\noindent{\bf Lemme #1. }{\it #2}\par}
\long\def\prop #1#2 {\par\noindent{\sc Proposition #1.--- }{\it
#2}\par}
\long\def\definition #1#2 {\par\noindent{\sc Définition
#1.--- }{\it #2}\par}
 \long\def\Th #1#2{\par\noindent{\sc Théorème
#1.--- }{\it #2}\par}

\long\def\Lemme#1{\par\noindent{\bf Lemme.--- }{\it
#1}\par}

\def\Rm{\par\noindent{{\it Remarque\,}: }}
\def\dem{\par\noindent{{\it Démonstration\,}: }}
\advance \parskip 4pt
\begin{document}
\newgeometry{left=4cm, right=4cm, top=2cm, bottom=3cm}
\title{Classification topologique des solutions du Problème
d'Apollonius}
\date{}
\author{Tchangang Tambekou Roger}
\renewcommand{\abstractname}{Abstract}
\maketitle {\abstract{We give a mathematical computation of the
number of solutions of Apollonius problem, by use of Lie Sphere
Geometry. Unlike in higher dimensions, the number of solutions
depends only on the topology of the configuration of the 3 objects.
It appears that our classification is non redundant, and far simpler
than those obtained previously } }

 Parmi les problèmes posés et résolus par Apollonius il y a
plus de 200 ans A.D, celui des contacts, qu'il aborde dans son Xème
traité est sans doute le plus célèbre. Il s'énonce comme suit :

\emph{Etant donnés 3 objets - cercles, points ou droites -
construire un objet tangent à tous les trois.}

\begin{figure}[http!]
 \begin{center}
    \includegraphics[scale=0.8]{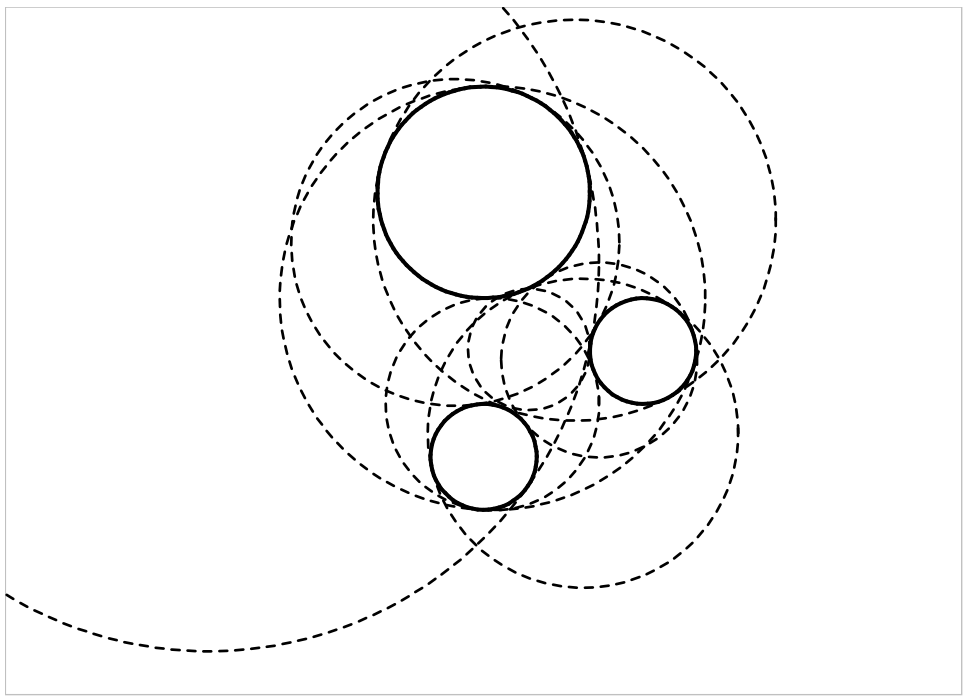}\\
\caption{Le nombre de solutions dépend de la configuration. Il y en
a 8 au maximum, comme sur cette figure.}\label{}

    \end{center}

        \end{figure}

    \restoregeometry
 La solution trouvée par Apollonius ayant disparu, à la fin du
XIVème, Viète \cite{u} s'est efforcé de la reconstituer.

Ce problème a fasciné plusieurs générations de mathématiciens,
depuis Pappus \cite{s}, Descartes \cite{g}, Newton \cite{r}, Euler
\cite{h}, Gauss \cite{l}, Cauchy \cite{c}, jusqu'à Coxeter \cite{f}.
Le problème d'Apollonius trouve des applications dans la théorie des
signaux, et dans plusieurs autres domaines de la mathématique
moderne \cite{e}. Des techniques géométriques ingénieuses, à la
règle et au compas \cite{k} ou bien par l'intermédiaire des
coniques, permettent de calculer et représenter les cercles
solutions. Mais ces outils n'assurent en rien que l'on a la totalité
des solutions, spécialement dans le cas où le nombre de solutions
n'est pas maximal.

Le problème d'Apollonius peut être généralisé à une dimension $n$
supérieure, en demandant que les objets trouvés soient tangents à
n+1 objets. Mais dans ce cas, le nombre de solutions ne dépend plus
seulement de la topologie des configurations d'objets, mais aussi de
leur équivalence au sens de la géométrie de Moebius.

Ce n'est qu'à la fin du XIXème que Muirhead \cite{q} pose le
problème - lui aussi intéressant - du nombre de solutions selon les
configurations des cercles, et donne sa propre classification des
configurations de 3 cercles, qui malheureusement est redondante et
incomplète. Un siècle plus tard, par des techniques de géométrie
inversive, Bruen, Fisher et Wilker \cite{b} en donnent une
simplification en 1983. En 1995, Fillmore et Paluzny \cite{i}
donnent la leur, en utilisant des quadriques projectives. Ces
classifications restent assez fastidieuses.

Enfin, R.D. Knight \cite{n} en 2005 donne une solution intéressante
au problème d'Apollonius orienté, c'est-à-dire au cas où les cercles
sont orientés et où l'on impose que les orientations des tangentes
soient les mêmes au point de contact. La technique de Knight utilise
la géométrie de Laguerre et celle de Moebius, qui sont des cas
particuliers de la géométrie de la Sphère de Lie (qui est le cadre
naturel d'un tel problème). Malheureusement, Knight ne relie pas le
problème dans le cas orienté au cas général, tel que posé par
Apollonius.

\widowpenalty10000

L'objet de cet article est de donner une technique simple pour
calculer le nombre de solutions pour une configuration donnée, à
partir de ses intersections et de la présence ou non de cercles
séparants et de points de tangence.

Nous utiliserons principalement l'inversion euclidienne et ses
propriétés, pour donner une classification  topologique agréable.
Nous introduirons les notions de puissance d'un cercle orienté par
rapport à un autre, et de discriminant de trois cercles orientés.

La technique utilisée ici est une spécialisation de raisonnements
plus généraux relevant du cadre de la géométrie de la Sphère de Lie,
même si elle n'est jamais citée.

Je remercie chaleureusement Théophile Noulaquape qui a vérifié les
calculs et réalisé les figures.

\section{La classification topologique}
Si une configuration de 3 objets au sens d'Apollonius contient des
droites, par une inversion, on la transforme en une configuration de
cercles ou de points. Nous ne considérerons donc que des
configurations avec des cercles de rayon nul ou pas, ce qui nous
évite d'introduire les notions d'intersection à l'infini et de
tangence à l'infini.

\subsection{Propriétés des configurations}
\subsubsection{Cercles Orientés}
 Dans le plan orienté, un cercle est dit non trivial s'il
est de rayon non nul. Un cercle définit un disque ouvert et un
disque fermé de même rayon, dont il est la frontière.

On désignera un cercle non orienté du plan par une lettre majuscule.
Le cercle $C$ admet 2 orientations. Les lettres minuscules
désigneront des cercles orientés ; Soit $x$ le centre du cercle $C$
et $r$ son rayon. Soit $c$ le cercle $C$ muni de l'orientation
positive, par exemple. On dira que $c$ est de centre $x$ et de rayon
$r$, et que $\overline{c}$ est de centre $x$ et de rayon $-r$. Ainsi
les cercles orientés négativement ont un rayon négatif.

L'application\footnote{cette application confère à l'espace des
cercles orientés une structure d'espace de Lorentz, sur lequel
agit le groupe affine de Lorentz.} qui à un cercle orienté $c$
associe $(x,r)$ est une bijection de l'ensemble des cercles
orientés vers l'espace $\mathbb{R}^3$.
\subsubsection{Cercle séparant}
\definition{1}{On dira que le cercle $C$ (ou $c$)
entoure le cercle $D$(ou $d$) si le disque ouvert défini par $C$ (ou
$c$) contient le disque ouvert défini par $D$ (ou $d$).

Lorsque le disque ouvert défini par $C$ contient le disque fermé
défini par $D$, on dit que le cercle $C$ entoure strictement le
cercle $D$.}

\definition{2}{Soit $(C_1, C_2, C_3)$ une configuration de 3 cercles du
plan euclidien, $i\in \{1,2,3\}$. Le cercle $C_i$ est dit
strictement séparant lorsque les deux autres cercles sont contenus
dans des composantes connexes distinctes du plan privé de $C_i$.

Le cercle $C_i$ est dit séparant (au sens large) lorsque les deux
autres cercles sont contenus dans les adhérences des composantes
connexes distinctes du plan privé de $C_i$.}

Le cas de cercles séparants sans être strictement séparants apparaît
seulement lorsqu'on a des cercles tangents.

\begin{figure}[http!]
 \begin{center}
    \includegraphics[scale=0.8]{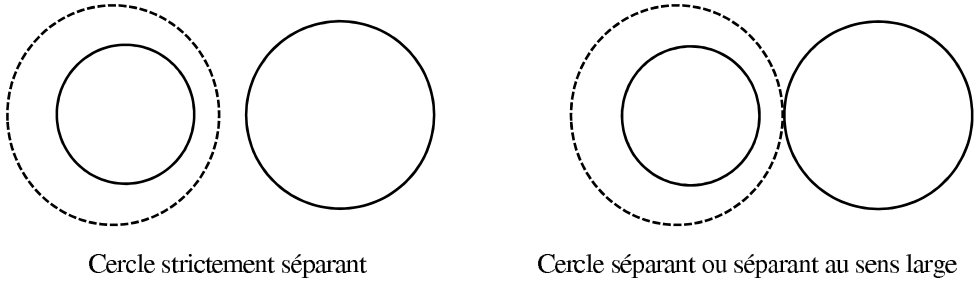}\\
\caption{Cercle séparant}\label{}
    \end{center}
    \end{figure}

\definition{3}{Un point d'intersection des cercles d'une configuration est dit double s'il appartient à tous les 3 cercles de la configuration.}
Une configuration de trois cercles contient au maximum un cercle
séparant.
\subsection{Classification topologique des solutions du Problème\\
d'Apollonius} Le tableau suivant donne la classification des
solutions du problème d'Apollonius, selon la topologie de la
configuration, c'est-à-dire son nombre de points d'intersection,
de points doubles, et la présence ou non de cercles séparants ou
de cercles tangents.

Il est évident que si une configuration contient un cercle
strictement séparant, alors il n'y a pas de solution au problème
d'Apollonius.

\begin{center}
\begin{tabular}{|p{7cm}||c|c|c|c|c|c|c|}
\hline \bf Nombre de solutions &\multicolumn{7}{c|}{\bf points d'intersection}\\
\hline
\bf La configuration contient &0&1&2&3&4&5&6\\
\hline
\hline
un cercle strictement séparant&0&&&&&&\\
\hline
pas de cercles tangents et pas de points doubles&8&&4&&4&&8\\
\hline un ou deux points de tangence, pas de cercle
séparant, pas de points doubles&&6&5&4&5&6&\\
\hline trois points de tangence&&&&5&&&\\
\hline des cercles tangents et un cercle
séparant, pas de points doubles&&2&3&&&&\\
\hline
des points doubles&&$\infty$&2&3&5&&\\
\hline
\end{tabular}
\end{center}

Les cases vides du tableau représentent des cas impossibles.

On voit qu'on a en tout 17 cas distincts et que
\begin{itemize}
    \item Le nombre maximal de solutions est 8.
    \item Il n'existe pas de configuration à 7 solutions.
    \item Il n'existe pas de configuration avec une solution
    unique.
\end{itemize}

Pour démontrer  cette classification, nous aurons recours à une
étude plus détaillée des cercles orientés. Nous résoudrons d'abord
un problème d'Apollonius orienté.

\subsection{Propriétés des cercles orientés}
\subsubsection{Tangence orientée}
 \definition{4}{ Deux cercles orientés $c$ et $d$ sont dits tangents
lorsque $C$ et $D$ sont tangents et qu'ils admettent la même
tangente orientée au point de contact.

Si le cercle $C$ est réduit à un point $x$, alors $c$ est tangent à
$d$ lorsque $x\in C$.}
\begin{figure}[http!]
 \begin{center}
    \includegraphics[scale=0.8]{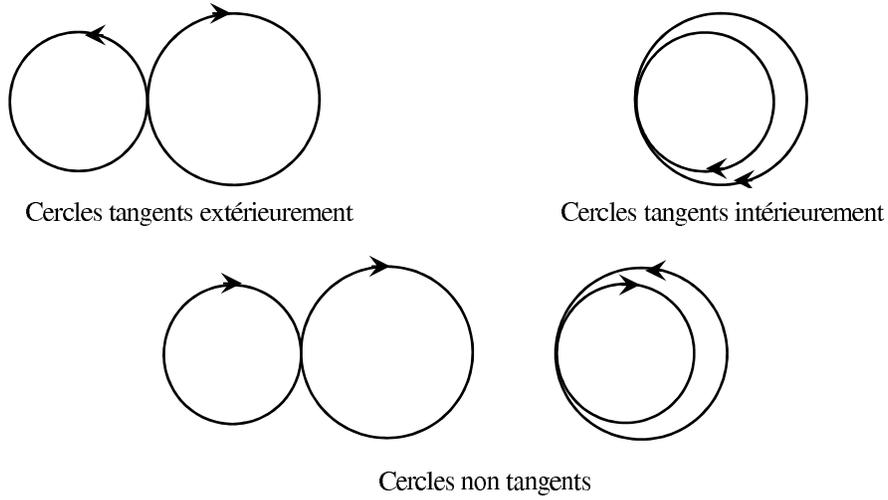}\\
\caption{ Cercles tangents, cercles non tangents}\label{}
    \end{center}
    \end{figure}
\newpage
\subsubsection{Puissance d'un cercle orienté par rapport à un
autre} Nous introduisons  cette notion qui généralise celle de
puissance d'un point par rapport à un cercle, en vue de définir le
discriminant de 3 cercles orientés. Comme dans le cas de la
puissance d'un point par rapport à un cercle, cette notion a un sens
géométrique précis.

\definition{5}{Soit $c_1(x_1,r_1)$ et $c_2(x_2,r_2)$ des cercles orientés. La
puissance de $c_1$ par rapport à $c_2$ est
$$P(c_1,c_2)=\|x_1-x_2\|^2 - (r_1-r_ 2)^2=Q(c_1-c_2), $$
où $Q(c)=\|x\|^2-r^2$, la forme quadratique de
Lorentz\footnote{Cette forme quadratique est invariante par le
groupe de Lorentz, Mais pas par le groupe affine de Lorentz.} de
signature $(2,1)$.}

\Rm lorsqu'aucun des 2 cercles n'entoure l'autre, $P(c_1,c_2)$ est
égal au carré de la distance entre $c_1$ et $c_2$, prise le long de
l'une des deux tangentes orientées communes.
\begin{figure}[http!]
 \begin{center}
    \includegraphics[scale=0.8]{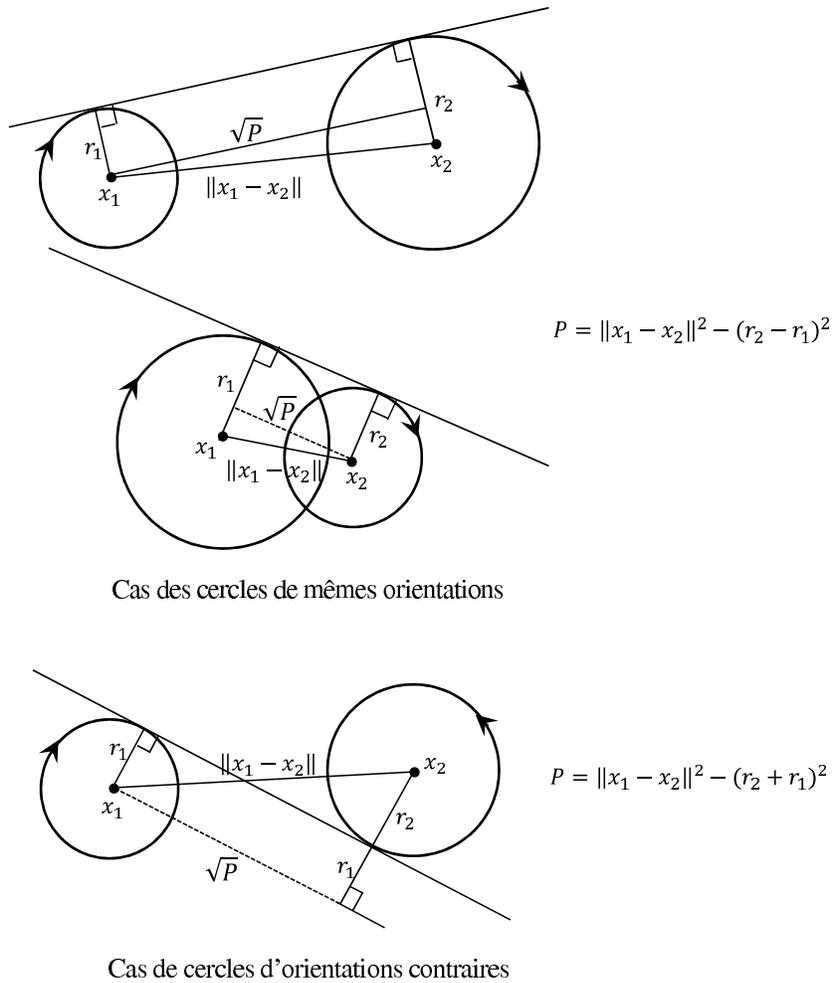}\\
\caption{Puissance de $c_1$ par rapport à $c_2$.}\label{}
    \end{center}
    \end{figure}
\newpage
La notion de puissance d'un cercle orienté par rapport à un autre
donne des informations précises sur la position relative des deux
cercles. Par exemple, si les deux cercles sont disjoints et qu'aucun
n'entoure l'autre,\\alors $P(c_1,c_2)> 0$, quelles que soient les
orientations.

\begin{figure}[http!]
 \begin{center}
    \includegraphics[scale=0.8]{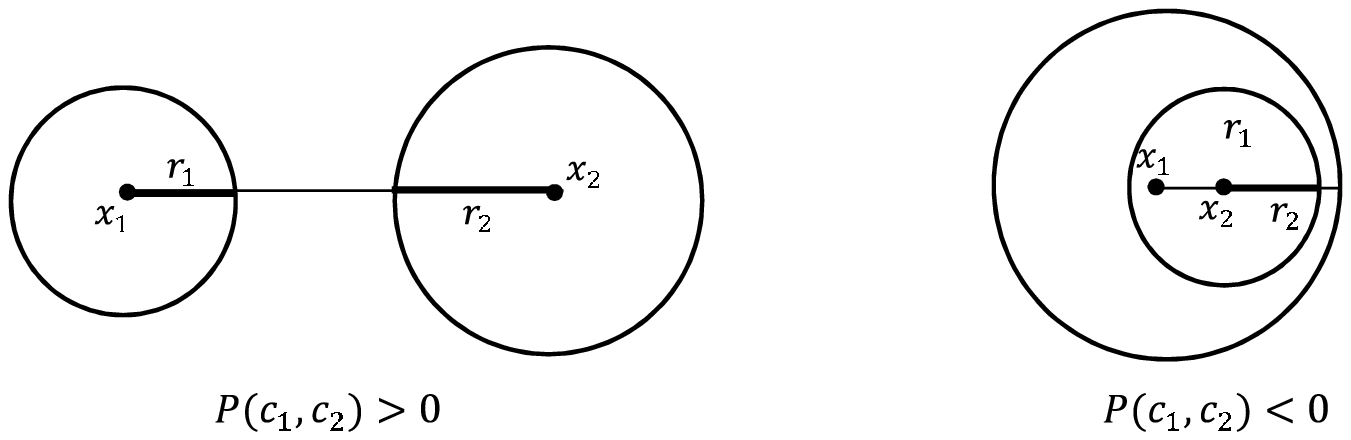}\\
    \end{center}
    \end{figure}

 Par contre,
si l'un des cercles entoure strictement l'autre, alors
$P(c_1,c_2)< 0$, indépendamment des orientations. Si les cercles
se coupent, on a deux cas

\begin{figure}[http!]
 \begin{center}
    \includegraphics[scale=0.8]{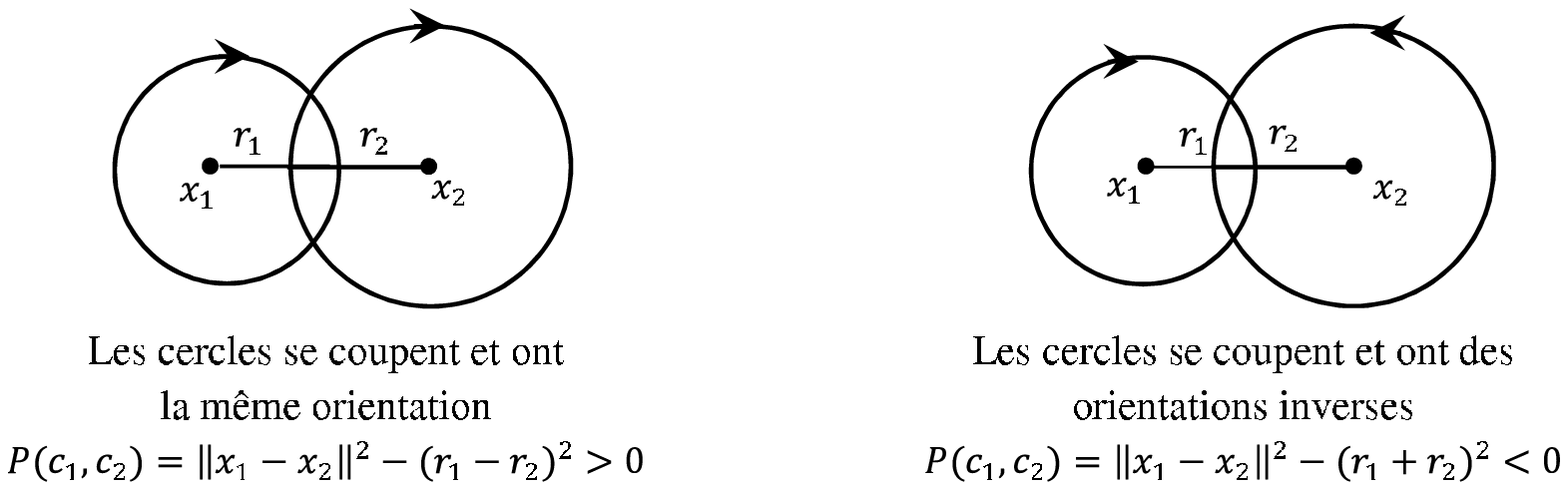}\\
    \end{center}
    \end{figure}

\prop{1}{Soient $c_1$ et $c_2$ deux cercles orientés du plan. Les
cercles $c_1$ et $c_2$ sont tangents si et seulement si
$P(c_1,c_2)= 0$, C'est-à-dire si et seulement si la puissance de
$c_1$ par rapport à $c_2$ est nulle.}

\dem On suppose les deux cercles tangents. Les centres $x_1$, $x_2$
et le point de tangence sont alignés, d'où le résultat. Inversement,
on considère la droite $(x_1x_2)$ qui est un axe de symétrie pour
les deux cercles. On la munit d'une graduation, en sorte que
$x_1\geq x_2$. On sait que $d(x_1, x_2)=\varepsilon(r_1-r_2)$, où
$\varepsilon=\pm1$. On a alors $x_1 -\varepsilon r_1=
x_2-\varepsilon r_2$, qui est le point de tangence.\par La position
relative de $c_1$ et $c_2$ est entièrement déterminée par le calcul
de $P(c_1,c_2)$ et $P(c_1, \overline{c}_2)$. Ces résultats sont
contenus dans le tableau suivant :

$$\begin{array}{|l|c|c|}
 \hline
P(c_1,c_2)&P(c_1, \overline{c}_2)&\mbox{Conclusion}\\
\hline
 +&+&\parbox[t]{8cm}{$c_1$ et
$c_2$ sont disjoints, aucun n'entourant l'autre} \\\hline
 +&-&\parbox[t]{8cm}{$c_1$ et $c_2$ se coupent et ont même sens} \\ \hline
 -&+& \parbox[t]{8cm}{$c_1$ et $c_2$ se coupent et ont des sens
contraires}\\ \hline
 -&-& \parbox[t]{8cm}{l'un des cercles entoure l'autre}\\ \hline
 0&+&\parbox[t]{8cm}{les deux cercles sont tangents extérieurement} \\ \hline
 0&-& \parbox[t]{8cm}{les deux cercles sont tangents intérieurement}\\ \hline
\end{array}$$

Ce tableau permet d'affirmer : \Th{1}{Les cercles orientés $c_1$ et
$c_2$ admettent exactement deux tangentes orientées communes si et
seulement si $P(c_1,c_2) > 0$.

Les cercles orientés $c_1$ et $c_2$ admettent exactement une
tangente orientée commune si et seulement si $P(c_1,c_2) = 0$.

Les cercles orientés $c_1$ et $c_2$ n'admettent aucune tangente
orientée commune si et seulement si $P(c_1,c_2) < 0$.}

\dem En effet, Lorsque les cercles orientés se coupent, les
tangentes communes n'existent que si les cercles ont la même
orientation.

\section{Problème d'Apollonius Orienté}
 Revenons aux 3 cercles non orientés
$C_1$, $C_2$ et $C_3$. Soit $D$ un cercle (ou une droite)  tangent à
ces 3 cercles. Si on munit $D$ d'une orientation, alors on peut
orienter les cercles $C_i$ en sorte que $d$ soit tangent aux $c_i$,
$i=1,2,3$. Donc $d$ ou $\overline{d}$ (et par suite $D$) est tangent
à l'un des quatre triplets de cercles orientés suivants : $(c_1,
c_2, c_3)$, $(c_1, c_2, \overline{c}_3)$, $(c_1, \overline{c}_2,
c_3)$, $(c_1, \overline{c}_2,\overline{c}_3)$. Soient $S_1, S_2,
S_3$ et $S_4$ les ensembles de solutions respectifs. Plus
précisément, $S_3$ est l'ensemble des cercles non orientés $D$ tels
que $d$ ou $\overline d$ est tangent à $c_1$, $\overline c_2$ et
$c_3$.

On rappelle que $d$ désigne le cercle $D$ muni d'une orientation.

Nous appellerons Problème d'Apollonius Orienté la recherche du
cardinal de l'un des ensembles $S_i$. Pour avoir toutes les
solutions du problème d'Apollonius, il suffit de considérer la
réunion des solutions aux Problèmes d'Apollonius Orientés obtenus en
prenant successivement toutes les orientations possibles des cercles
de la configuration.

\Lemme {Si les cercles $c_i$ sont de rayon non nul et n'ont aucun
point double (i.e. $C_1\cap C_2\cap C_3=\emptyset$), alors les
ensembles $S_i$, $i=1,2,3,4$, sont disjoints deux à deux.

Si les 3 cercles ont un ou deux points communs, alors ces points
appartiennent à tous les ensembles $S_i$. }

\dem Soit $c$ un cercle orienté appartenant à $S_1\cap S_2$, par
exemple. Alors $c$ est tangent à $c_3$ et à $\overline{c}_3$, et par
suite, $c$ est de rayon nul, ce qui est impossible, car les trois
cercles $c_i$ n'ont pas de point commun.

\definition{6}{Discriminant de trois cercles}
Soient $c_1, c_2$ et $c_3$ des cercles orientés du plan euclidien.
On suppose qu'ils ne sont pas tous les trois tangents en un même
point. On appelle discriminant de ces trois cercles orientés le
nombre
$$D=Q(c_1-c_2)Q(c_1-c_3)Q(c_2-c_3).$$
\Th {2}{Le Problème d'Apollonius Orienté a exactement 2 solutions si
$D
> 0$; une solution unique si $D=0$ et il n'admet pas de solution
si $D < 0$.}

Pour démontrer ce théorème, on va d'abord se ramener au cas où l'un
des cercles est réduit à un point, et ensuite calculer le signe du
discriminant selon les configurations possibles. Soit $j$ un vecteur
de $\mathbb{R}^3$. Les cercles orientés $c$ et $d$ sont tangents si
et seulement si les cercles orientés $c-j$ et $d-j$ le sont. En
effet, On a $Q((c-j)-(d-j))=Q(c-d)=0$. Par suite, si on considère
$j=(0,0,r_1)$, le cercle orienté $c_1-j$ est réduit à un point, et
on peut donc, sans nuire à la généralité, considérer que $c_1$ est
réduit à un point. C'est d'ailleurs la technique utilisée par
plusieurs auteurs pour construire une solution à la règle et au
compas, à la différence que plusieurs n'abordent pas la question des
rayons négatifs. La démonstration ci-dessous s'adapte facilement au
cas où $c_2$ ou $c_3$ sont de rayons nuls.

\dem Considérons une inversion $I$ de centre $x_1$, donc qui envoie
le point $c_1=x_1$ à l'infini. Si le discriminant est non nul, alors
aucun des trois cercles n'est tangent à un autre. Donc le point
$x_1$ n'appartient ni à $c_2$, ni à $c_3$. Par conséquent, $I(c_2)$
et $I(c_3)$ sont des cercles. Les solutions du problème
d'Apollonius, quand elles existent, sont transformées par $I$ en des
droites, tangentes communes à $I(c_2)$ et $I(c_3)$. Nous allons donc
considérer les différentes configurations possibles des 3 cercles,
calculer dans chaque cas le discriminant, et le nombre $N$ de
tangentes qu'on obtient en appliquant la transformation $I$, à
partir du signe de $P(I(c_2),I(c_3))$. Il faut remarquer qu'une
inversion de centre $x_1$ renverse l'orientation des cercles ne
contenant pas $x_1$, mais préserve celle des cercles contenant
$x_1$, ce qui permet de construire les tableaux suivants.

\begin{figure}[http!]
 \begin{center}
    \includegraphics[scale=0.8]{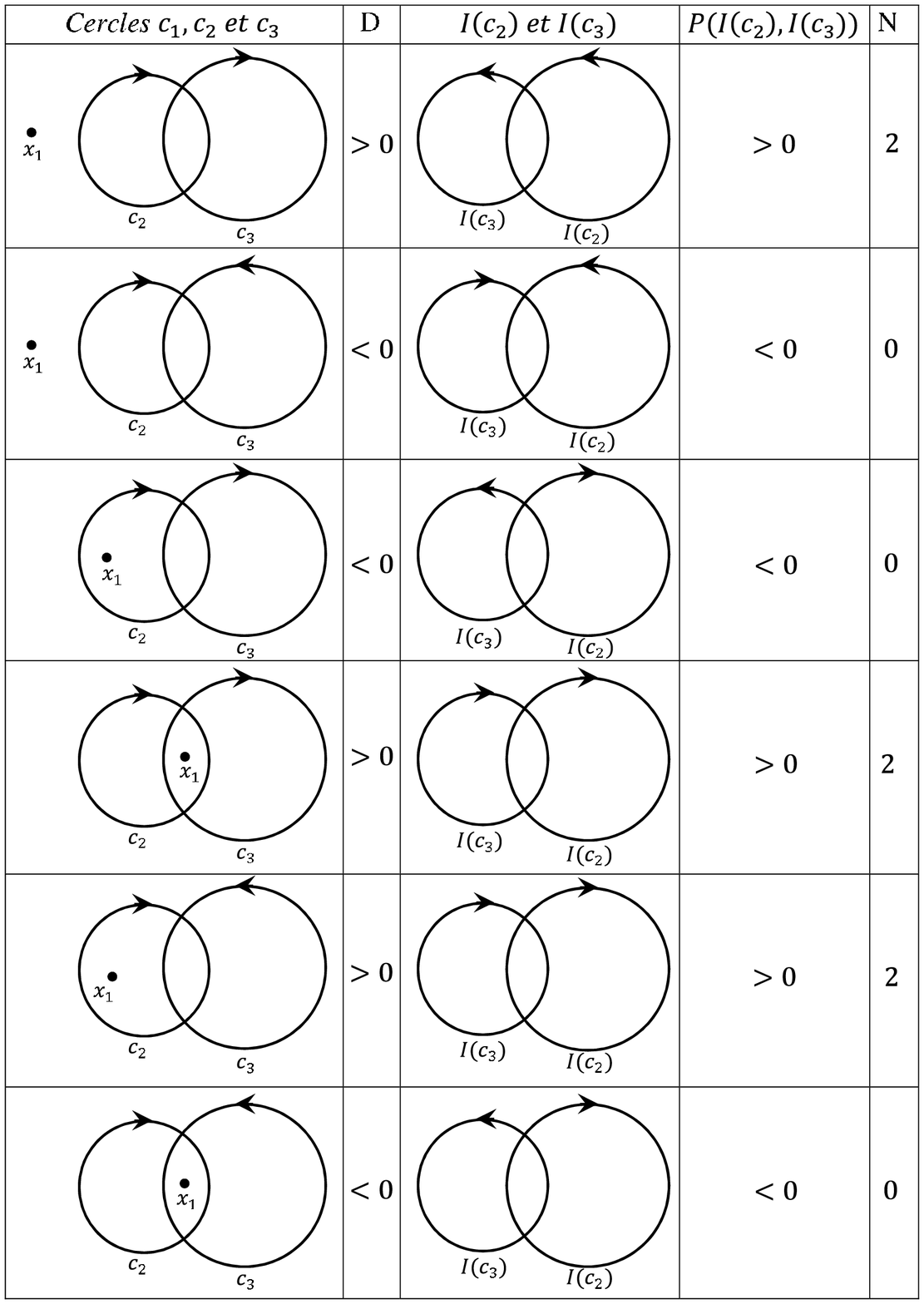}\\
\caption{Cas où les cercles $c_2$ et $c_3$ se coupent}\label{}
    \end{center}
    \end{figure}

\begin{figure}[http!]
 \begin{center}
    \includegraphics[scale=0.8]{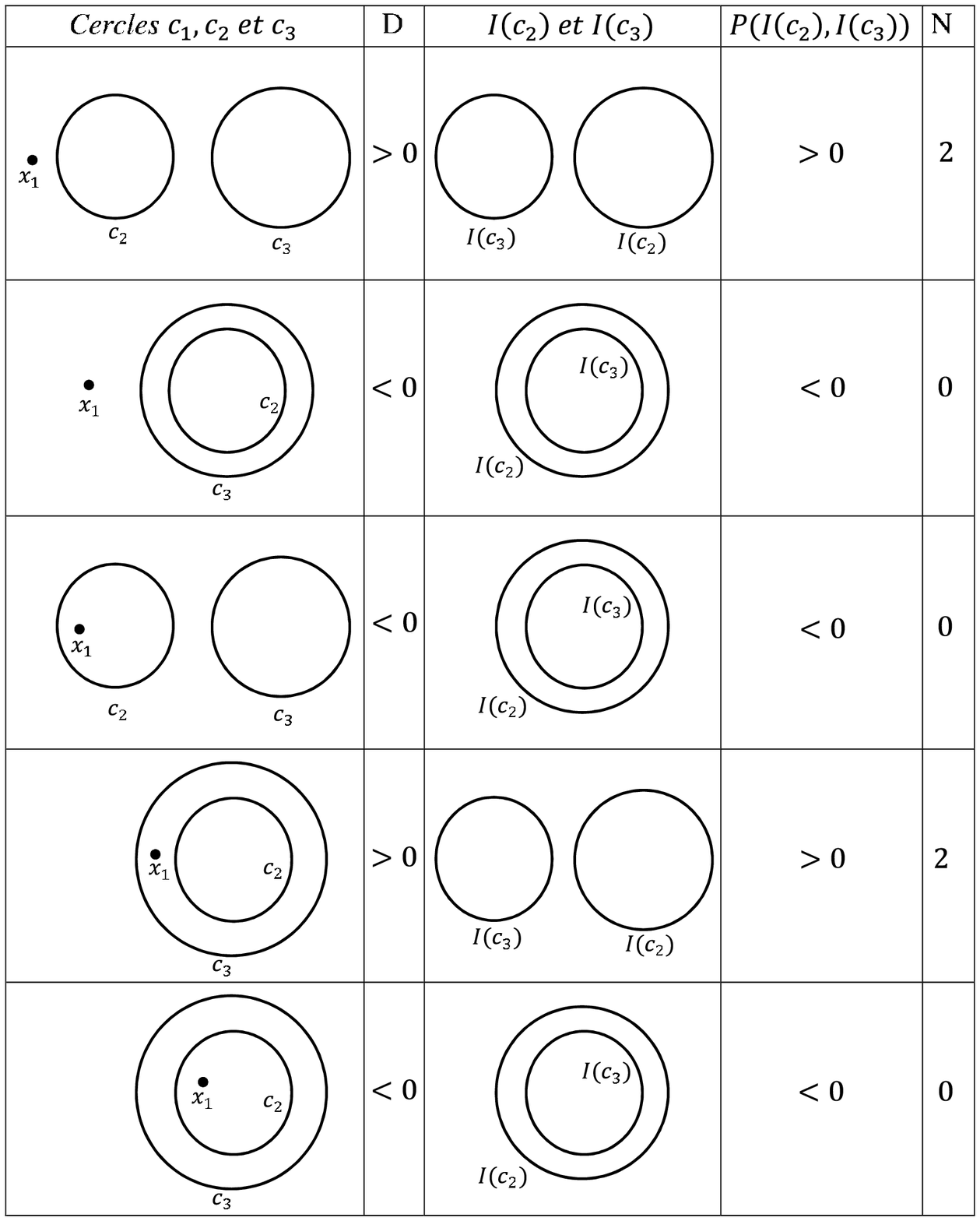}\\
\caption{Cas où les cercles $c_2$ et $c_3$ sont disjoints}\label{}
    \end{center}
    \end{figure}

\newpage
La section suivante est consacrée au cas où le discriminant est nul
(deux des cercles orientés au moins sont tangents).
 \subsection{Cas de cercles orientés tangents}
\subsubsection{Deux des cercles orientés seulement sont tangents}
On appelle $c_1$ le cercle orienté qui n'est tangent à aucun autre.
Les cercles $I(c_2)$ et $I(c_3)$ sont tangents, et admettent une et
une seule tangente commune qui passe par le point de contact. Le
Problème d'Apollonius Orienté a une seule solution dans ce cas.
\subsubsection{$c_1$ est tangent à $c_2$ qui est tangent à $c_3$}
Dans ce cas, le point $c_1$ appartient à l'un des cercles, mettons
$c_2$. L'inversion $I$ transforme le cercle $c_2$ en une droite
$\Delta$. Tout cercle tangent à $c_2$ et passant par $x_1$ est
transformé par l'involution $I$ en une droite parallèle à $\Delta$.
Comme il existe une seule telle tangente orientée à $I(c_3)$, on en
déduit que le Problème d'Apollonius a une solution unique.

\begin{figure}[http!]
 \begin{center}
    \includegraphics[scale=0.8]{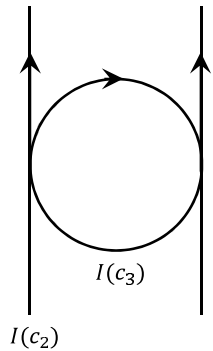}\\
\caption{Tangente commune}\label{}
    \end{center}
    \end{figure}

\subsubsection{Les 3 cercles orientés sont tangents deux à deux}
Dans ce cas, ils ne peuvent pas être tous tangents extérieurement,
car ils seraient tous d'orientations différentes, ce qui est
impossible. Donc si deux des cercles sont tangents extérieurement,
le troisième est tangent intérieurement à l'un d'eux, donc il y a un
seul point d'intersection pour les 3 cercles; de même, si les trois
cercles sont tous tangents intérieurement, alors il y a un seul
point d'intersection. Tout cercle passant par le point commun et
tangent à l'un des cercles est solution du problème d'Apollonius, et
dans les deux cas, le Problème d'Apollonius admet une infinité de
solutions.\par
\begin{figure}[http!]
 \begin{center}
    \includegraphics[scale=0.8]{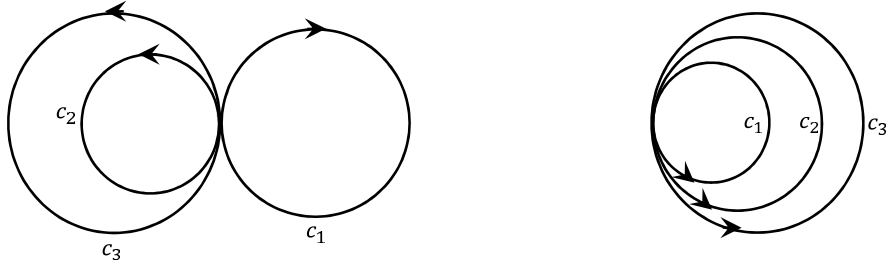}\\
\caption{Les deux configurations sont homéomorphes via une inversion
par rapport à $c_3$.}\label{}
    \end{center}
    \end{figure}
Donc lorsque deux des cercles orientés sont tangents (discriminant
nul) et qu'il n'y a pas de point double, on a une seule solution au
Problème d'Apollonius Orienté, ce qui démontre le théorème.

\section{Nombre de solutions du Problème d'Apollonius}
Nous pouvons maintenant compter les solutions du problème
d'Apollonius, par une technique de type do it yourself. Soient
$c_1$, $c_2$ et $c_3$ les cercles orientés positivement associés à
$C_1$, $C_2$ et $C_3$ respectivement. \\Soit par exemple la
configuration suivante : $C_1$ tangent à $C_2$ tangent à $C_3$
tangent à $C_1$, où les cercles non orientés sont tous tangents,
avec des points de tangence distincts. On calcule les cardinaux des
ensembles $S_1$, $S_2$, $S_3$ et $S_4$. On a le tableau suivant :

\begin{center}
\begin{tabular}{|l|c|c|c|c|}
\hline
Configuration&Discriminant&N\\
\hline $c_1$,$c_2$,$c_3$&+&2\\ \hline $c_1$,$\overline{c}_2$,$c_3$&0&1\\
\hline $c_1$,$c_2$,$\overline{c}_3$&0&1\\ \hline $c_1$,$\overline{c}_2$,$\overline{c}_3$&0&1\\
\hline $C_1$,$C_2$,$C_3$&TOTAL&5\\ \hline
\end{tabular}
\end{center}

\begin{figure}[http!]
 \begin{center}
    \includegraphics[scale=0.8]{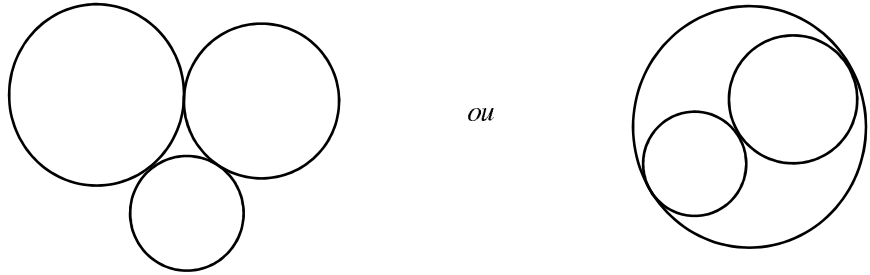}\\
\caption{On passe d'une figure à l'autre par une inversion.}\label{}
    \end{center}
    \end{figure}
\newpage
Les cercles sont tangents 2 à 2, les deux configurations se
déduisent l'une de l'autre par une inversion, et donnent le même
nombre de cercles tangents. Remarquons que chacun des trois cercles
est lui-même solution.

\subsection{Classification}
\subsubsection{L'un des cercles est strictement séparant}
 Alors, quelles que soient les orientations, on a $D < 0$,
et par suite, comme on s'y attend, il n'y a pas de solution.

\subsubsection{Les
cercles sont disjoints sans cercle séparant}
 Soit aucun cercle n'entoure un autre, soit  l'un des cercles entoure les deux autres
; Dans ce cas $D > 0$, quelles que soient les orientations des
cercles, et on a 8 solutions.

\subsubsection{Deux des cercles $C_2$ et $C_3$
seulement se coupent}
 Le troisième cercle $C_1$ est disjoint des deux
autres. Remarquons que  le produit $Q(c_1-c_2)Q(c_1-c_3)$ garde le
même signe quelles que soient les orientations des cercles. On a
donc le tableau suivant à un changement de signe près, qui donne 4
solutions :

\begin{center}
\begin{tabular}{|l|c|c|c|c|}
\hline
Configuration&Discriminant&N\\
\hline $c_1$,$c_2$,$c_3$&+&2\\ \hline $c_1$,$\overline{c}_2$,$c_3$&-&0\\
\hline $c_1$,$c_2$,$\overline{c}_3$&-&0\\ \hline $c_1$,$\overline{c}_2$,$\overline{c}_3$&+&2\\
\hline $C_1$,$C_2$,$C_3$&TOTAL&4\\ \hline
\end{tabular}
\end{center}

\subsubsection{Le cercle $C_1$ coupe $C_2$ et $C_2$ coupe $C_3$}
 Les cercles $C_1$
et $C_3$ sont disjoints, et l'un entoure l'autre, ou non. Dans les
deux cas, on obtient le tableau suivant à un changement de signe du
discriminant près :

\begin{center}
\begin{tabular}{|l|c|c|c|c|}
\hline
Configuration&Discriminant&N\\
\hline $c_1$,$c_2$,$c_3$&+&2\\ \hline $c_1$,$\overline{c}_2$,$c_3$&+&2\\
\hline $c_1$,$c_2$,$\overline{c}_3$&-&0\\ \hline $c_1$,$\overline{c}_2$,$\overline{c}_3$&-&0\\
\hline $C_1$,$C_2$,$C_3$&TOTAL&4\\ \hline
\end{tabular}
\end{center}

\subsubsection{Tous les cercles se coupent deux à deux, en des points
distincts}
 On
a le tableau suivant :

\begin{center}
\begin{tabular}{|l|c|c|c|c|}
\hline
Configuration&Discriminant&N\\
\hline $c_1$,$c_2$,$c_3$&+&2\\ \hline $c_1$,$\overline{c}_2$,$c_3$&+&2\\
\hline $c_1$,$c_2$,$\overline{c}_3$&+&2\\ \hline $c_1$,$\overline{c}_2$,$\overline{c}_3$&+&2\\
\hline $C_1$,$C_2$,$C_3$&TOTAL&8\\ \hline
\end{tabular}
\end{center}

\subsubsection{Les trois cercles ont un point commun, sans être
tangents} Alors le point double appartient aux 4 ensembles de
solutions, et par suite, le nombre de solutions est $8-4+1 = 5$.
C'est le cas par exemple de 3 droites qui forment un triangle dans
le plan, quand on ajoute le point à l'infini à cette configuration
qu'on transforme par une inversion.

\subsubsection{Les trois cercles ont deux points communs}
Ces deux points doubles appartiennent tous aux 4 ensembles  de
solutions et par suite le nombre de solutions est $8-8+2=2$.

\subsubsection{Deux des cercles sont tangents, le
troisième ($c_1$) est disjoint des deux autres} Si le troisième
cercle entoure les deux autres ou pas, le produit
$Q(c_1-c_2)Q(c_1-c_3)$ garde un signe constant quelles que soient
les orientations des cercles. On obtient à homéomorphisme près deux
possibilités selon l'existence ou non d'un cercle séparant.

\begin{center}
\begin{tabular}{|l|c|c|c|c|}
\hline
Configuration&Discriminant&N\\
\hline $c_1$,$c_2$,$c_3$&+&2\\ \hline $c_1$,$\overline{c}_2$,$c_3$&0&1\\
\hline $c_1$,$c_2$,$\overline{c}_3$&0&1\\ \hline $c_1$,$\overline{c}_2$,$\overline{c}_3$&+&2\\
\hline $C_1$,$C_2$,$C_3$&TOTAL&6\\ \hline
\end{tabular}
\quad
\begin{tabular}{|l|c|c|c|c|}
\hline
Configuration&Discriminant&N\\
\hline $c_1$,$c_2$,$c_3$&-&0\\ \hline $c_1$,$\overline{c}_2$,$c_3$&0&1\\
\hline $c_1$,$c_2$,$\overline{c}_3$&0&1\\ \hline $c_1$,$\overline{c}_2$,$\overline{c}_3$&-&0\\
\hline $C_1$,$C_2$,$C_3$&TOTAL&2\\ \hline
\end{tabular}
\end{center}
\begin{figure}[http!]
 \begin{center}
    \includegraphics[scale=0.8]{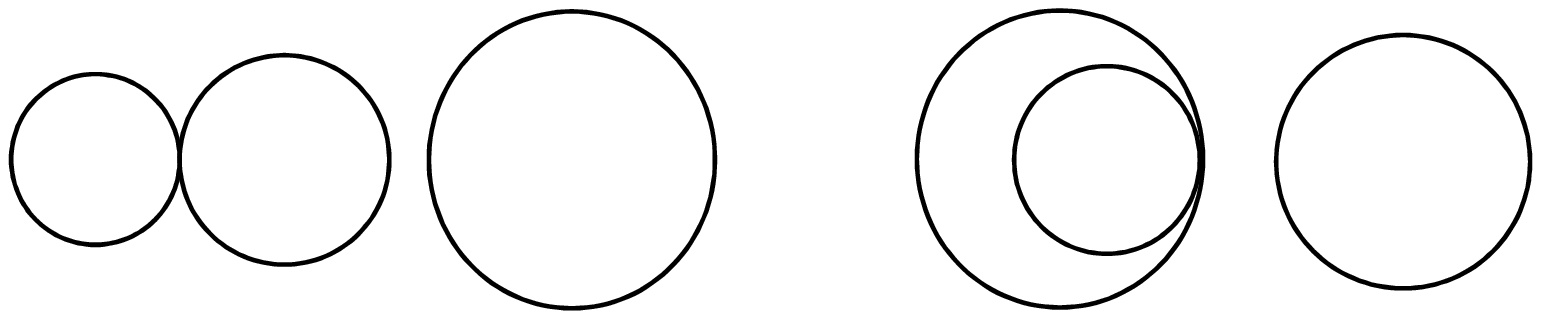}\\
\caption{}\label{}
    \end{center}
    \end{figure}
\vspace{-0.5cm}
\subsubsection{Deux des cercles sont tangents, et le troisième coupe
l'un d'eux seulement} Ici encore, on a deux situations qui se
déduisent l'une de l'autre par une inversion. On a le tableau
suivant :

\begin{center}
\begin{tabular}{|l|c|c|c|c|}
\hline
Configuration&Discriminant&N\\
\hline $c_1$,$c_2$,$c_3$&-&0\\ \hline $c_1$,$\overline{c}_2$,$c_3$&0&1\\
\hline $c_1$,$c_2$,$\overline{c}_3$&0&1\\ \hline $c_1$,$\overline{c}_2$,$\overline{c}_3$&+&2\\
\hline $C_1$,$C_2$,$C_3$&TOTAL&4\\ \hline
\end{tabular}
\end{center}

\subsubsection{Deux des cercles sont tangents, et le troisième
($C_1$) coupe les deux autres en dehors du point de contact}

\begin{center}
\begin{tabular}{|l|c|c|c|c|}
\hline
Configuration&Discriminant&N\\
\hline $c_1$,$c_2$,$c_3$&+&2\\ \hline $c_1$,$\overline{c}_2$,$c_3$&0&1\\
\hline $c_1$,$c_2$,$\overline{c}_3$&0&1\\ \hline $c_1$,$\overline{c}_2$,$\overline{c}_3$&+&2\\
\hline $C_1$,$C_2$,$C_3$&TOTAL&6\\ \hline
\end{tabular}
\end{center}

\subsubsection{Deux des cercles sont tangents, et le troisième coupe
les 2 autres au point d'intersection}
 Dans ce cas, le point d'intersection
apparaît chaque fois comme tangent aux quatre configurations. Ce
qui donne $6-4+1=3$ cercles tangents, le point d'intersection
inclus.

\subsubsection{Le cercle $C_1$ est tangent à $C_2$ qui est tangent à
$C_3$} Ici, on a une chaîne de 3 cercles tangents. On obtient  l'un
des 2 tableaux suivants selon que la configuration a un cercle
séparant ou non.

\begin{center}
\begin{tabular}{|l|c|c|c|c|}
\hline
Configuration&Discriminant&N\\
\hline $c_1$,$c_2$,$c_3$&+&2\\ \hline $c_1$,$\overline{c}_2$,$c_3$&0&1\\
\hline $c_1$,$c_2$,$\overline{c}_3$&0&1\\ \hline $c_1$,$\overline{c}_2$,$\overline{c}_3$&0&1\\
\hline $C_1$,$C_2$,$C_3$&TOTAL&5\\ \hline
\end{tabular}
\quad
\begin{tabular}{|l|c|c|c|c|}
\hline
Configuration&Discriminant&N\\
\hline $c_1$,$c_2$,$c_3$&0&1\\ \hline $c_1$,$\overline{c}_2$,$c_3$&0&1\\
\hline $c_1$,$c_2$,$\overline{c}_3$&-&0\\ \hline $c_1$,$\overline{c}_2$,$\overline{c}_3$&0&1\\
\hline $C_1$,$C_2$,$C_3$&TOTAL&3\\ \hline
\end{tabular}
\end{center}
\begin{figure}[http!]
 \begin{center}
    \includegraphics[scale=0.8]{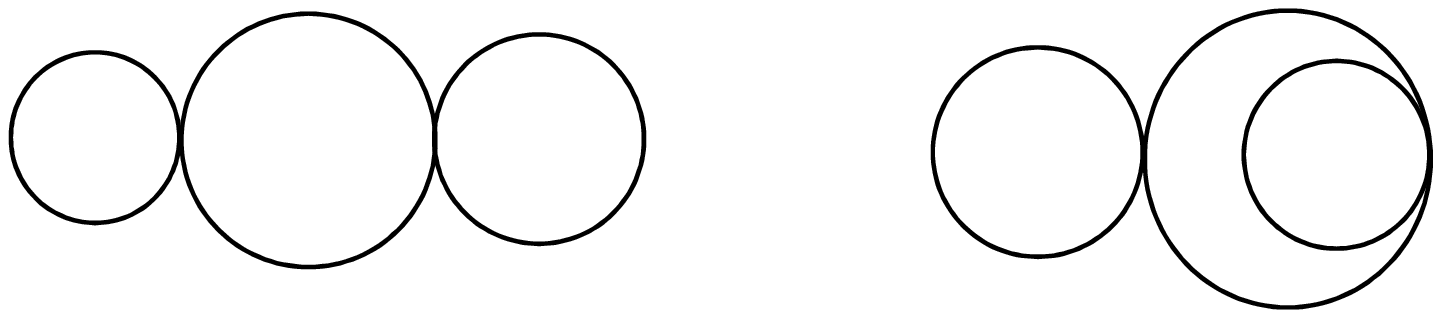}\\
\caption{}\label{}
    \end{center}
    \end{figure}
Le résultat ne change pas lorsque $C_3$ coupe $C_1$.
\subsubsection{Les 3 cercles sont tangents 2 à 2} Ce cas a déjà été
traité plus haut, et on a trouvé 5 solutions.

\subsection{Cas d'un point et deux cercles}
La méthode ici consiste à envoyer le point à l'infini et à compter
ensuite les tangentes communes aux deux cercles.

\begin{center}
\begin{tabular}{|p{8cm}||c|c|c|c|}
\hline \bf Nombre de solutions &\multicolumn{4}{c|}{\bf points d'intersection}\\
\hline
\bf La configuration contient &0&1&2&3\\
\hline \hline
un cercle strictement séparant&0&&&\\
\hline
pas de cercles tangents ni séparants et pas de points doubles&4&2&2&2\\
\hline deux cercles tangents, pas de cercle
séparant, pas de points doubles&&3&2&\\
\hline un points de tangence et un cercle séparant&&1&&\\
\hline
des points doubles&&$\infty$&1&\\
\hline
\end{tabular}
\end{center}

\subsection{Cas de 2 points et d'un cercle}
La méthode est la même que celle utilisée dans le cas précédent et
on obtient les résultats suivants:
\begin{center}
\begin{tabular}{|p{8cm}||c|}
\hline  \bf Configuration &\bf Nombre de solutions\\
\hline \hline
Le cercle sépare strictement les deux points&0\\
\hline
un  des points au moins est tangent au cercle&1\\
\hline le cercle ne sépare pas les points et aucun point n'appartient au cercle&2\\
\hline
\end{tabular}
\end{center}

\subsection{Cas de 3 points}
Il existe une solution unique, le cercle circonscrit au triangle
défini par les trois points ou la droite passant par ces points.

\vskip 2mm \scriptsize
\parindent -8pt
\begin{tabular}{l@{\qquad\qquad\qquad\qquad}l}
\it Address : \rm Tchangang Tambekou Roger&\it  Teaching address :\\
\verb"roger_tchangang@yahoo.fr"&Département de Mathématiques \\
Centre de Recherches en Géométrie et Applications&Faculté des Sciences\\
B.P. 8 451 Yaoundé&Université de Yaoundé I\\
Cameroun&Cameroun
\end{tabular}
\end{document}